\documentclass[reqno,12pt]{amsart}
\usepackage{amsmath,amsthm,amssymb,amsfonts,amscd,graphicx}
\input xy
\xyoption{all}
\theoremstyle{plain} 
	\newtheorem{thm}{Theorem}[section]
	\newtheorem*{thm*}{Theorem}
	
	\newtheorem{lem}[thm]{Lemma}
	
	\newtheorem{prop}[thm]{Proposition}
	\newtheorem*{prop*}{Proposition}
	
	\newtheorem*{conj*}{Conjecture}
	
\theoremstyle{definition}
	\newtheorem{defn}[thm]{Definition}

\theoremstyle{remark}
	\newtheorem{rem}[thm]{Remark}
	
	\newtheorem*{pf}{Proof}
\numberwithin{equation}{section}

\def\RR{{\mathbb R}}

\def\ZZ{{\mathbb Z}}

\def\I{{\mathcal I}}

\def\R{{\mathcal R}}

\begin{document}
\title{On Simply-Laced Generalized Root Systems}
\date{\today}
\author{Shunsuke Nakamura}
\address{Department of Mathematics, Graduate School of Science, Osaka University, 
Toyonaka Osaka, 560-0043, Japan}
\email{u891956f@gmail.com}
\author{Yuuki Shiraishi}
\address{Department of Mathematics, Graduate School of Science, Kyoto University, 
Kyoto, 606-8202, Japan}
\email{yshiraishi@math.kyoto-u.ac.jp}
\author{Atsushi Takahashi}
\address{Department of Mathematics, Graduate School of Science, Osaka University, 
Toyonaka Osaka, 560-0043, Japan}
\email{takahashi@math.sci.osaka-u.ac.jp}
\begin{abstract}
We show the uniqueness and existence of the Euler form for a simply-laced generalized root system. 
This enables us to show that the Coxeter element for a simply-laced generalized root system is
admissible in the sense of R.~W.~Carter. 
As an application, the isomorphism classes of simply-laced generalized root systems with 
positive definite Cartan forms are classified by Cartar's admissible diagrams associated to their Coxeter elements.
\end{abstract}
\maketitle
\maketitle
\section{Introduction}
The notion of simply-laced generalized root systems is introduced by Kyoji Saito \cite{saito1998}
in his study of period mappings of primitive forms in the deformation theory of isolated
hypersurface singularities. This is an axiomatization of the tuple consisting of
the middle homology group of the Milnor fibration, the intersection form, the set of vanishing cycles and 
the Milnor monodromy (see Definition~\ref{object}). 
One of the most famous examples of them is 
from the deformation theory of a simple elliptic singularity, which is called an elliptic root system.
The elliptic root systems and the Lie algebras associated to them are beyond the theory of Kac--Moody Lie algebras.
This class is extensively studied by many researchers including K. Saito.
We suggest only \cite{saito,st:1} here for the relevance to the present paper.  
A finite root system, an affine root systems, or the set of real roots in the lattice with 
the Cartan form of a Kac--Moody Lie algebra, which is called a KM root system here, also appears
as a substructure of a simply-laced generalized root system.
The difference between simply-laced generalized root systems and such root systems
is the Coxeter element, which plays an important role.

Certain simply-laced generalized systems are constructed from the deformation theory of affine cusp polynomials 
$x_{1}+x_{2}^{a_{2}}+x_{3}^{a_{3}}-cx_{1}x_{2}x_{3}$ (see \cite{ist:3, stw, T}). 
If we forget their Coxeter elements, then they all give the affine root system of type $\widetilde{A}_{a_{2}+a_{3}}$.
These simply-laced generalized root systems have the different Coxeter elements corresponding to different pairs $(a_{2},a_{3})$
with $a_{2}\le a_{3}$ and hence they are non-isomorphic to each other.
Therefore, simply-laced generalized root systems have finer information than the usual ones.
Due to the Coxeter elements, simply-laced generalized root systems can be equipped with  
not only representation theoretic information but also geometric one.

Except for studies of elliptic root systems mentioned above, general properties 
of simply-laced generalized root systems are hardly known.
However, it is shown by \cite{stw} that they are constructed from triangulated 
categories with good algebraic origins; the tuple consisting of the Grothendieck group,  
the symmetrized Euler form, the set of isomorphic classes of exceptional collections and 
the automorphism induced by the Serre functor.  
It is important to note here that the Euler form, the definitely important pairing, gives the Cartan form. 
Therefore, it is natural to ask whether a simply-laced generalized root system has such a nice bilinear form.

The main result of the present paper is the following positive answer to this$:$
\begin{thm}[Theorem~\ref{Euler form}]
Let $R=(K_0(R),I_R,\Delta_{re}(R),c_R)$ be a simply-laced generalized root system.
There exists a unique non-degenerate bilinear form $\chi_R:K_0(R)\times K_0(R)\longrightarrow\ZZ$ 
satisfying the following properties$:$
\begin{enumerate}
\item
We have $I_R=\chi_R+\chi_R^T$, namely,  
\begin{subequations}
\begin{equation*}
I_R(\lambda,\lambda')=\chi_R(\lambda,\lambda')+\chi_R(\lambda',\lambda),\quad \lambda,\lambda'\in K_0(R).
\end{equation*}
\item We have 
\begin{equation*}
\chi_R(\lambda,\lambda')=-\chi_R(\lambda',c_R\lambda),\quad \lambda,\lambda'\in K_0(R).
\end{equation*}
In particular, the Coxeter transformation $c_R$ is equal to $-\chi_R^{-1}\chi_R^{T}$.
\end{subequations}
\end{enumerate}
\end{thm}
If the Cartan form of a simply-laced generalized root system is positive definite,
we can show that the Coxeter element is admissible (see Definition~\ref{w admi}) 
in the sense of R. W. Carter (see Proposition~\ref{admissible}) by Theorem~\ref{Euler form}.
Here we define morphisms between simply-laced generalized root systems as  
isometric $\ZZ$--homomorphisms which are compatible with Coxeter transformations and
induce maps between the sets of real roots (see Definition~\ref{morphism}).
Then we can classify isomorphism classes of simply-laced generalized root systems
whose Cartan forms are positive definite by admissible diagrams introduced and classified by R. W. Carter$:$
\begin{thm}[Theorem~\ref{classification thm}]
Let $\R^{{\rm irr, pd}}$ be the set of irreducible positive definite generalized root systems,
$\Gamma$ the set of diagrams in Proposition~\ref{classifycation} {\rm (ii)} 
and $A_{R}$ the diagram in Proposition~\ref{classifycation} {\rm (ii)} which is 
isomorphic to the admissible diagram associated to the Coxeter element $c_{R}$ of $R$.
The map
\begin{equation*}
\varphi:\R^{{\rm irr, pd}}\longrightarrow \Gamma, \quad R\mapsto A_{R},
\end{equation*}
induces the  natural bijection$:$
\begin{equation*}
\overline{\varphi}:\R^{{\rm irr, pd}}/\sim \overset{\cong}\longrightarrow \Gamma, \quad [R]\mapsto A_{R},
\end{equation*}
where $[R]$ is the isomorphism class of $R \in \R^{{\rm irr, pd}}$.
\end{thm}

\bigskip
\noindent
{\it Acknowledgment}\\
\indent
The first and second named authors would like to thank Akishi Ikeda and Tadashi Ishibe
for valuable discussions and encouragement.
The second named author is supported by Research Fellowships of Japan Society for 
the Promotion for Young Scientists.
The third named author is supported by JSPS KAKENHI Grant Number 24684005.
\section{Simply-laced generalized root systems}

\subsection{Simply-laced generalized root systems and morphisms}
We recall the definition of simply-laced generalized root systems
introduced by K. Saito \cite{saito1998}. In addition, we define a certain class of morphisms between
simply-laced generalized root systems.

\begin{defn}\label{object}
A {\it simply-laced generalized root system $R$} of rank $\mu$ consists of  
\begin{itemize}
\item
a free $\ZZ$-module $K_0(R)$ of rank $\mu$ called the {\it root lattice},
\item
a symmetric bilinear form $I_R:K_0(R)\times K_0(R)\longrightarrow\ZZ$ called the {\it Cartan form},
\item
a subset $\Delta_{re}(R)$ of $K_0(R)$ called the {\it set of real roots},
\item
an element $c_R$ of $W(R)$ called the {\it Coxeter transformation}, 
\end{itemize}
such that
\begin{enumerate}
\item
$K_0(R)=\ZZ\Delta_{re}(R)$.
\item
For all $\alpha\in\Delta_{re}(R)$, $I_{R}(\alpha,\alpha)=2$.
\item 
For all $\alpha\in\Delta_{re}(R)$, the element $r_\alpha$ of ${\rm Aut}(K_0(R),I_R)$, 
the group of automorphisms of $K_0(R)$ respecting $I_R$, defined by
\begin{equation*}
r_\alpha(\lambda):=\lambda-I_R(\lambda,\alpha)\alpha,\quad \lambda\in K_0(R)
\end{equation*}
makes $\Delta_{re}(R)$ invariant, namely, $r_\alpha(\Delta_{re}(R))=\Delta_{re}(R)$.
\item 
Let $W(R)$ be the {\it Weyl group} of $R$ defined by
\begin{equation*}
W(R):=\langle r_{\alpha}~\vert~\alpha\in\Delta_{re}(R)\rangle\subset {\rm Aut}(K_0(R),I_R)\rangle.
\end{equation*}
Then there exists an ordered set $B=(\alpha_1,\dots, \alpha_\mu)$ of elements of $\Delta_{re}(R)$ called a {\it root basis} of $R$ 
which satisfies 
\begin{enumerate}
\item
$K_0(R)=\displaystyle\bigoplus_{i=1}^\mu \ZZ\alpha_i$, 
\item
$W(R)=\langle r_{\alpha_1},\dots,r_{\alpha_\mu}\rangle$, 
\item
$\Delta_{re}(R)=W(R)B$,
\item
$c_R=r_{\alpha_1}\cdots r_{\alpha_\mu}$.
\end{enumerate}
\end{enumerate}
\end{defn}

\begin{defn}
Let $R=(K_0(R),I_R,\Delta_{re}(R),c_R)$ be a simply-laced generalized root system. 
$R$ is called {\it reducible} if there exist subsets $\Delta_1, \Delta_2\subset \Delta_{re}(R)$
satisfying the following conditions$:$
\begin{equation*}
\begin{split}
&\Delta_{re}(R)=\Delta_1\cup \Delta_2,\ \Delta_1\cap \Delta_2=\emptyset, \\
&I_R (\alpha_1, \alpha_2)=0 \quad \text{for} \quad \alpha_1\in\Delta_1, \ \alpha_2\in\Delta_2.
\end{split}
\end{equation*} 
A simply-laced generalized root system $R$ is called {\it irreducible} if $R$ is not reducible. 
\end{defn}

\begin{defn}
Put 
\begin{equation*}
{\rm rad}(I_{R}):=\left\{\lambda\in K_0(R)~\vert~I_{R}(\lambda,\lambda')=0 \quad \text{ for all }\lambda'\in K_0(R)\right\}
\end{equation*}
and call it the {\rm radical} of $R$.
\end{defn}
\begin{lem}\label{lem2.4}
Let $R=(K_0(R),I_R,\Delta_{re}(R),c_R)$ be a simply-laced generalized root system.
\begin{enumerate}
\item
If $\alpha\in\Delta_{re}(R)$, then $-\alpha\in\Delta_{re}(R)$.
\item
Suppose that $\alpha\in\Delta_{re}(R)$. Then $c\alpha\in\Delta_{re}(R)$ for some constant $c\in\ZZ$ if and only if $c=\pm 1$.
\item
If $\alpha\in\Delta_{re}(R)$, then $r_{-\alpha}=r_{\alpha}$.
\end{enumerate}
\end{lem}
\begin{pf}
It is almost clear from the facts that $-\alpha=r_{\alpha}(\alpha)$, $I_R(c\alpha,c\alpha)=2c^2$ and the definition of reflections.
\qed
\end{pf}
\begin{prop}
Let $R=(K_0(R),I_R,\Delta_{re}(R),c_R)$ be a simply-laced generalized root system.
Let $I:K_0(R)\times K_0(R)\longrightarrow\ZZ$ be a symmetric bilinear form which satisfies the following conditions$:$
\begin{itemize}
\item
For all $\alpha\in\Delta_{re}(R)$, we have $I(\alpha,\alpha)=2$.
\item
For all $\alpha\in\Delta_{re}(R)$, we have $I(r_{\alpha}(\lambda),r_{\alpha}(\lambda'))=I(\lambda,\lambda')$ 
for all $\lambda,\lambda'\in K_0(R)$.
\end{itemize}
Then we have $I=I_R$.
\end{prop}
\begin{pf}
For all $\lambda\in K_0(R)$ and $\alpha\in\Delta_{re}(R)$, we have  
\begin{eqnarray*}
I(\lambda,\alpha)
&=&I(r_{\alpha}(\lambda),r_{\alpha}(\alpha))=-I(r_{\alpha}(\lambda),\alpha)\\
&=&-I\left(\lambda-I_R(\lambda,\alpha)\alpha,\alpha\right)=-I(\lambda,\alpha)+2 I_R(\lambda,\alpha).
\end{eqnarray*}
Therefore, we have $I(\lambda,\alpha)=I_R(\lambda,\alpha)$ and hence $I=I_R$.
\qed
\end{pf}

\begin{defn}\label{morphism}
Let $R=(K_0(R),I_R,\Delta_{re}(R),c_R)$ and $R'=(K_0(R'),I_{R'},\Delta_{re}(R'),c_{R'})$ 
be simply-laced generalized root systems.
A {\it morphism} $\phi:R\longrightarrow R'$ is a $\ZZ$-homomorphism $\phi:K_0(R)\longrightarrow K_0(R')$ satisfying the following 
conditions$:$
\begin{enumerate}
\item
For all $\lambda_1,\lambda_2\in K_0(R)$, we have 
\begin{equation*}
I_{R'}\left(\phi(\lambda_1),\phi(\lambda_2)\right)=I_R(\lambda_1,\lambda_2).
\end{equation*}
\item
For all $\alpha\in \Delta_{re}(R)$, we have $\phi(\alpha)\in\Delta_{re}(R')$.
\item
We have $\phi\circ c_{R}=c_{R'}\circ \phi$, namely, 
The following diagram commutes$:$
\begin{equation*}
\xymatrix{
K_0(R) \ar[r]^{\phi}\ar[d]_{c_R} &K_0(R')\ar[d]^{c_{R'}}\\
K_0(R) \ar[r]^{\phi} &K_0(R')
}.
\end{equation*}
\end{enumerate}
\end{defn}
\begin{rem}
The notion of isomorphisms of simply-laced generalized root systems is defined in the obvious way. 
Namely, a morphism $\phi:R\longrightarrow R'$ is an isomorphism if there exists a morphism $\psi:R'\longrightarrow R$ 
such that $\psi\circ \phi={\rm id}_R$ and $\phi\circ \psi={\rm id}_{R'}$.
\end{rem}

\begin{lem}
Let $R=(K_0(R),I_R,\Delta_{re}(R),c_R)$ and $R'=(K_0(R'),I_{R'},\Delta_{re}(R'),c_{R'})$ 
be simply-laced generalized root systems and $\phi:R\longrightarrow R'$ a morphism.
If $\lambda\in K_0(R)$ satisfies $\phi(\lambda)=0$, then $\lambda$ belongs to ${\rm rad}(I_R)$.
\end{lem}
\begin{pf}
For all $\lambda'\in K_0(R)$, we have 
\[
I_R(\lambda,\lambda')=I_{R'}\left(\phi(\lambda),\phi(\lambda')\right)=0.
\]
Therefore this lemma holds.
\qed
\end{pf}
\begin{lem}\label{commute c and phi}
Let $R=(K_0(R),I_R,\Delta_{re}(R),c_R)$ and $R'=(K_0(R'),I_{R'},\Delta_{re}(R'),c_{R'})$ 
be simply-laced generalized root systems and $\phi:R\longrightarrow R'$ a morphism.
For $\alpha\in\Delta_{re}(R)$, we have 
\begin{equation*}
\phi \circ r_{\alpha}=r_{\phi(\alpha)}\circ \phi.
\end{equation*}
\end{lem}
\begin{pf}
For any $\lambda\in K_0(R)$, we have 
\begin{eqnarray*}
\phi\left(r_{\alpha}(\lambda)\right)&=&\phi\left(\lambda-I_R(\lambda,\alpha)\alpha\right)\\
&=&\phi(\lambda)-I_R(\lambda,\alpha)\phi(\alpha)\\
&=&\phi(\lambda)-I_{R'}\left(\phi(\lambda),\phi(\alpha)\right)\phi(\alpha)\\
&=&r_{\phi(\alpha)}\left(\phi(\lambda)\right).
\end{eqnarray*}
Therefore this lemma holds.
\qed
\end{pf}
\begin{lem}
Let $R=(K_0(R),I_R,\Delta_{re}(R),c_R)$ and $R'=(K_0(R'),I_{R'},\Delta_{re}(R'),c_{R'})$ 
be simply-laced generalized root systems and $\phi:R\longrightarrow R'$ a morphism.
The $\ZZ$-submodule ${\rm Im}(\phi)\subset K_0(R')$ is invariant under the Coxeter transformation $c_{R'}:$
\begin{equation*}
c_{R'}\left({\rm Im}(\phi)\right)\subset {\rm Im}(\phi).
\end{equation*}
\end{lem}
\begin{pf}
It is obvious from the definition $c_{R'}\circ \phi =\phi\circ c_R$.
\qed
\end{pf}

\subsection{Euler form} In this subsection, we show the existence and uniqueness of the Euler form
on a simply-laced generalized root system.

\begin{thm}\label{Euler form}
Let $R=(K_0(R),I_R,\Delta_{re}(R),c_R)$ be a simply-laced generalized root system.
There exists a unique non-degenerate bilinear form $\chi_R:K_0(R)\times K_0(R)\longrightarrow\ZZ$ 
satisfying the following properties$:$
\begin{enumerate}
\item
We have $I_R=\chi_R+\chi_R^T$, namely,  
\begin{subequations}
\begin{equation}\label{2.4a}
I_R(\lambda,\lambda')=\chi_R(\lambda,\lambda')+\chi_R(\lambda',\lambda),\quad \lambda,\lambda'\in K_0(R).
\end{equation}
\item We have 
\begin{equation}\label{2.4b}
\chi_R(\lambda,\lambda')=-\chi_R(\lambda',c_R\lambda),\quad \lambda,\lambda'\in K_0(R).
\end{equation}
In particular, the Coxeter transformation $c_R$ is equal to $-\chi_R^{-1}\chi_R^{T}$.
\end{subequations}
\end{enumerate}
\end{thm}
\begin{pf}
First, we show the existence of such a bilinear form. 
We prove it by exactly the same argument in Proposition 2.4 in \cite{H and K}$:$
\begin{lem}\label{exist of Euler}
Let $(\alpha_1, \dots, \alpha_{\mu})$ be a root basis of a simply-laced generalized root system $R$. 
Define the bilinear form 
$\chi_{R}$ on $K_{0}(R)$ by the upper triangular part of the Cartan matrix $I_{R}$, namely,
\begin{equation*}
\chi_{R}(\alpha_{i},\alpha_{j}):=
\begin{cases}
I_{R}(\alpha_{i},\alpha_{j}) \quad \text{if} \quad i<j,\\
1 \quad \text{if} \quad i=j,\\
0 \quad \text{if} \quad i>j.
\end{cases}
\end{equation*}
Then the bilinear form $\chi_{R}$ satisfies the properties in Theorem~\ref{Euler form}. 
\end{lem}
\begin{pf}
Set $E_{i}:=(\alpha_{1},\dots,\alpha_{i})$ for $i=1,\dots, \mu$. 
We show this lemma by the induction for $i$.
For $\displaystyle y=m \alpha_{1} \in \ZZ\alpha_{1}$, we have
\begin{eqnarray*}
\chi_{R}(y, r_{\alpha_{1}}(x))
&=&\chi_{R}(y,x)-I_{R}(x,\alpha_{1})\cdot \chi_{R}(y,\alpha_{1})\\
&=&\chi_{R}(y,x)-I_{R}(x,\alpha_{1})\cdot m\\
&=&\chi_{R}(y,x)-I_{R}(x,y)=-\chi_{R}(x,y).
\end{eqnarray*}

Assume that $\chi_{R}(y, r_{\alpha_{1}}\cdots r_{\alpha_{i}}(x))=-\chi_{R}(x, y)$ for $\displaystyle y \in \bigoplus^{i}_{j=1}\ZZ\alpha_{j}$.
Under this assumption, we shall show that 
$\chi_{R}(y, r_{\alpha_{1}}\cdots r_{\alpha_{i+1}}(x))=-\chi_{R}(x, y)$ for $\displaystyle y \in \bigoplus^{i+1}_{j=1}\ZZ\alpha_{j}$.
In order to show this, it is sufficient to consider the following two cases$;$
$\displaystyle y \in \bigoplus^{i}_{j=1}\ZZ\alpha_{j}$ and $y \in \ZZ\alpha_{i+1}$. 
For the case that $\displaystyle y \in \bigoplus^{i}_{j=1}\ZZ\alpha_{j}$, we have
\begin{eqnarray*}
\chi_{R}(y, r_{\alpha_{1}}\cdots r_{\alpha_{i}}\cdot r_{\alpha_{i+1}}(x))
&=&-\chi_{R}(r_{\alpha_{i+1}}(x),y)\\
&=&-\chi_{R}(r_{\alpha_{i+1}}(x)-x,y)-\chi_{R}(x,y)\\
&=&-\chi_{R}(x,y),
\end{eqnarray*}
since $\displaystyle r_{\alpha_{i+1}}\cdots r_{\alpha_{\mu}}(x)-x\in \bigoplus^{\mu}_{j=i+1}\ZZ\alpha_{j}$.
For the case that $\displaystyle y=m \alpha_{i+1} \in \ZZ\alpha_{i+1}$, we have
\begin{eqnarray*}
\chi_{R}(y, r_{\alpha_{1}}\cdots r_{\alpha_{i}}\cdot r_{\alpha_{i+1}}(x))
&=&\chi_{R}(y, r_{\alpha_{1}}\cdots r_{\alpha_{i}}\cdot r_{\alpha_{i+1}}(x)-r_{\alpha_{i+1}}(x))+\chi_{R}(y, r_{\alpha_{i+1}}(x))\\
&=&\chi_{R}(y, r_{\alpha_{i+1}}(x))\\
&=&\chi_{R}(y,x)-I_{R}(x,\alpha_{i+1})\cdot \chi_{R}(y,\alpha_{i+1})\\
&=&\chi_{R}(y,x)-I_{R}(x,\alpha_{i+1})\cdot m\\
&=&\chi_{R}(y,x)-I_{R}(x,y)=-\chi_{R}(x,y)
\end{eqnarray*}
since $\displaystyle r_{\alpha_{1}}\cdots r_{\alpha_{i}}\cdot r_{\alpha_{i+1}}(x)-r_{\alpha_{i+1}}(x)\in \bigoplus^{i}_{j=1}\ZZ\alpha_{j}$.
Therefore this lemma holds.
\qed
\end{pf}

We show the following lemma necessary to prove the uniqueness.
\begin{lem}\label{gluing-1}
There exist  finitely many sub $\ZZ$--modules $V_{\bf i}$ indexed by a certain finite set $\I$
satisfying the following conditions$:$
\begin{enumerate}
\item $\displaystyle \sum_{{\bf i}\in \I} V_{\bf i}=K_{0}(R)$.
\item $V_{\bf i}\oplus {\rm rad}(I_{R})=K_{0}(R)$ for all ${\bf i}\in \I$.
In particular, the restriction of $I_{R}$ to $V_{\bf i}\times V_{\bf i}$ is non-degenerate for all ${\bf i}\in \I$. 
\end{enumerate}
\end{lem}
\begin{pf}
Set $\mu_{0}:={\rm rank} \ {\rm rad} (I_{R})$.
We shall construct the sub $\ZZ$--modules $V_{\bf i}$ explicitly.
It is obvious that there exists a $\ZZ$--basis $\beta_{1},\dots, \beta_{\mu}$ of $K_{0}(R)$ such that the restriction of $I_{R}$ to 
$\displaystyle \sum_{i=1}^{\mu-\mu_{0}}\ZZ\beta_{i}\times \sum_{i=1}^{\mu-\mu_{0}}\ZZ\beta_{i}$ is non-degenerate
and $I_{R}(\beta_{j},\beta_{j})\ne 0$ for all $j=1,\dots, \mu$. Indeed we can take $\beta_{1},\dots, \beta_{\mu}$
as an appropriate permutation of the initial root basis $(\alpha_{1},\dots,\alpha_{\mu})$.

A $\ZZ$--basis of ${\rm rad} \ (I_{R})$ can be written as follows$:$
\begin{equation*}
\delta_{i}:=\sum_{k=1}^{\mu} m_{ik} \beta_{k}, \quad i=1,\dots \mu_{0}.
\end{equation*}
If there exists a nontrivial sequence of integers $(k_{1},\dots k_{\mu_{0}})$ such that
\begin{equation}\label{indep-of-sum}
\sum^{\mu_{0}}_{j=1}k_{j}\left(\sum^{\mu}_{k=\mu-\mu_{0}+1}m_{jk}\beta_{k}\right)=0,
\end{equation}
then the equation \eqref{indep-of-sum} implies 
$\displaystyle \sum^{\mu_{0}}_{j=1}k_{j}\delta_{j}\in \sum_{i=1}^{\mu-\mu_{0}}\ZZ\beta_{i}$ by the equations 
\begin{equation}\label{before-transf}
\displaystyle \sum_{k=\mu-\mu_{0}+1}^{\mu} m_{jk} \beta_{k}=
\delta_{j}-\sum_{k=1}^{\mu-\mu_{0}} m_{jk} \beta_{k}, \quad j=1,\dots, \mu_{0}.
\end{equation}
The equation \eqref{indep-of-sum} is contradictory to the assumption for $\displaystyle \sum_{i=1}^{\mu-\mu_{0}}\ZZ\beta_{i}$.
Hence the elements $\displaystyle \sum^{\mu}_{k=\mu-\mu_{0}+1}m_{jk}\beta_{k}$ for $j=1, \dots, \mu_{0}$ are 
linearly independent.
Moreover, by the equation \eqref{before-transf}, we have the following equation:
\begin{equation}\label{generator} 
K_{0}(R)=\sum_{i=1}^{\mu-\mu_{0}}\ZZ\beta_{i}+\sum_{i=1}^{\mu_{0}}\ZZ\delta_{i}
=\sum_{i=1}^{\mu-\mu_{0}}\ZZ\beta_{i}+\sum_{i=1}^{\mu_{0}}\ZZ\left(\sum_{k=\mu-\mu_{0}+1}^{\mu} m_{ik} \beta_{k}\right).
\end{equation}
Thus the elements $\beta_{i}$ for $i=1,\dots, \mu-\mu_{0}$ and 
$\displaystyle \sum^{\mu}_{k=\mu-\mu_{0}+1}m_{jk}\beta_{k}$ for $j=1, \dots, \mu_{0}$
are generators of $K_{0}(R)$.

We consider the following $\ZZ$--homomorphism $f:K_{0}(R)\rightarrow K_{0}(R)$:
\begin{equation*}
\begin{cases}
f(\beta_{i})=\beta_{i} \quad \text{if} \quad i=1,\dots,\mu-\mu_{0},\\
f(\beta_{i+\mu-\mu_{0}})=\displaystyle\sum_{k=\mu-\mu_{0}+1}^{\mu} m_{ik} \beta_{k} \quad \text{if} \quad i=1,\dots,\mu_{0}.
\end{cases}
\end{equation*}
By the equation \eqref{generator}, we have ${\rm Coker}(f)=0$. 
We denote the representation matrix of $f$ restricted to the sub $\ZZ$--module
$\displaystyle \sum^{\mu}_{k=\mu-\mu_{0}}\ZZ\beta_{k}$ as follows$:$
\begin{equation*}
M':=
\begin{pmatrix}
a_{\mu-\mu_{0}+1 1} & \cdots & a_{\mu-\mu_{0}+1 j} & \cdots & a_{\mu-\mu_{0}+1 \mu_{0}}\\
\vdots & \ddots &        &        & \vdots \\
a_{k 1} &        & a_{k j} &        & a_{k \mu_{0}} \\
\vdots &        &        & \ddots & \vdots \\
a_{\mu 1} & \cdots & a_{\mu j} & \cdots & a_{\mu \mu_{0}}
\end{pmatrix}.
\end{equation*}
If the absolute value of ${\rm det}(M')$ is greater than one, the cokernel of $f$ is nontrivial and
has torsions by the theory of elementary divisors.  
  
By the above argument, the matrix $M'$ is an element in ${\rm GL}(\mu_{0},\ZZ)$.
Then we have the following equation:
\begin{equation*}\label{elementary-transf}
(\beta_{\mu-\mu_{0}+1},\dots, \beta_{\mu})=(\delta_{1},\dots \delta_{\mu_{0}})(M')^{-1}+(-\sum_{k=1}^{\mu-\mu_{0}} m_{1k} \beta_{k},\dots, 
-\sum_{k=1}^{\mu-\mu_{0}} m_{\mu_{0}k} \beta_{k})(M')^{-1}.
\end{equation*}
We set $L^{(1)}_{i}(\delta_{1}, \dots, \delta_{\mu_{0}})$ and $L^{(2)}_{i}(\beta_{1}, \dots, \beta_{\mu-\mu_{0}})$ as follows$:$
\begin{eqnarray*}
(L^{(1)}_{1}(\delta_{1}, \dots, \delta_{\mu_{0}}),\dots, L^{(1)}_{\mu_{0}}(\delta_{1}, \dots, \delta_{\mu_{0}}))
&:=&(\delta_{1},\dots \delta_{\mu_{0}})(M')^{-1}\\
(L^{(2)}_{1}(\beta_{1}, \dots, \beta_{\mu-\mu_{0}}),\dots,L^{(2)}_{\mu_{0}}(\beta_{1}, \dots, \beta_{\mu-\mu_{0}})
&:=&\displaystyle(-\sum_{k=1}^{\mu-\mu_{0}} m_{1k} \beta_{k},\dots, -\sum_{k=1}^{\mu-\mu_{0}} m_{\mu_{0}k} \beta_{k})(M')^{-1}.
\end{eqnarray*}
It is clear that the elements $L^{(1)}_{i}(\delta_{1}, \dots, \delta_{\mu_{0}})$ for $i=1,\dots, \mu_{0}$ form 
a $\ZZ$--basis of ${\rm rad}(I_{R})$ by their constructions. 

Since $\beta_{1},\dots,\beta_{\mu}$ are linearly independent, we can define the following set$:$
\begin{equation*}
\displaystyle {\rm Supp}\left(\sum_{k=1}^{\mu} m_{k} \beta_{k}\right):=\{i\in \{1,\dots,\mu\}|m_{i}\ne 0\}.
\end{equation*}
Since $I_{R}(\beta_{\mu-\mu_{0}+i},\beta_{\mu-\mu_{0}+i})\ne 0$,
we have ${\rm Supp}(L^{(2)}_{i}(\beta_{1}, \dots, \beta_{\mu-\mu_{0}}))\ne \emptyset$.
For $i=1,\dots, \mu_{0}$, we take numbers $p_{i}\in {\rm Supp}(L^{(2)}_{i}(\beta_{1}, \dots, \beta_{\mu-\mu_{0}}))$.
We check that the following sub $\ZZ$--modules $V_{(i,p_{i})}$ and
the sub $\ZZ$--module $V_{1}:=\displaystyle \sum_{i=1}^{\mu-\mu_{0}}\ZZ\beta_{i}$ 
are the desired ones: 
\begin{equation*}
\displaystyle V_{(i,p_{i})}:=\ZZ\beta_{\mu-\mu_{0}+i}+\sum_{1\le j\le \mu-\mu_{0}, j\ne p_{i}} 
\ZZ\beta_{j}, \quad i=1,\dots \mu_{0}.
\end{equation*}
A nonzero element in $V_{(i,p_{i})}$ and $L^{(1)}_{j}(\delta_{1}, \dots, \delta_{\mu_{0}})$ 
for all $j=1,\dots,\mu_{0}$ are linearly independent. 
Hence the restriction $I_{R}$ to $V_{(i,p_{i})}\times V_{(i,p_{i})}$ is non-degenerate. 
In addition, we have $\displaystyle \sum_{{\bf i}\in \I} V_{\bf i}=K_{0}(R)$  
and $V_{\bf i}\oplus {\rm rad}(I_{R})=K_{0}(R)$ for each element ${\bf i}\in \I$ where we set
$\I:=\{1, (1,p_{1}),\dots, (\mu_{0},p_{\mu_{0}})\}$.
\qed
\end{pf}

Next, we show the uniqueness of such a bilinear form. 
\begin{lem}\label{uniqueness-of-euler}
Assume that there exists a non-degenerate bilinear form $\chi$ satisfying \eqref{2.4a} and \eqref{2.4b}.
Then the bilinear form $\chi$ is unique.
\end{lem}
\begin{pf}
Set $B_{{\bf i},1}:=\{\beta_{j}\}_{\beta_{j}\in V_{\bf i}}$ and 
$B_{{\bf i}, 2}:=\{L^{(1)}_{j}(\delta_{1}, \dots$, $\delta_{\mu_{0}})\}_{1\le j\le \mu_{0}}$ 
as in Lemma~\ref{gluing-1}
and denote by $I_{\bf i}$, $X_{\bf i}$ and $C_{\bf i}$ the representation matrices of $I_{R}$, $\chi$ and $c_{R}$
concerning the $\ZZ$--basis $B_{{\bf i},1}\cup B_{{\bf i},2}$ respectively.

The bilinear form $\chi$ (resp. $I_{R}$) defines an element in ${\rm Hom}_{\ZZ}(K_{0}(R), {\rm Hom}_{\ZZ}(K_{0}(R),\ZZ))$ by
$a\mapsto (b\mapsto \chi(a, b))$ (resp. $a\mapsto (b\mapsto I_{R}(a, b))$).
Then $X^{-1}_{\bf i}$ is the representation matrix between the $\ZZ$--bases 
$B'_{{\bf i},1}\cup B'_{{\bf i},2}$ and $B_{{\bf i},1}\cup B_{{\bf i},2}$
where $B'_{{\bf i},1}$ and $B'_{{\bf i},2}$ are the dual $\ZZ$--basis of $B_{{\bf i},1}$ and $B_{{\bf i},2}$ respectively.

Denote by $(M)_{B,B'}$ the submatrix of the representation matrix $M$ corresponding to sub $\ZZ$--bases $B$ and $B'$. 
The submatrix $(I_{\bf i})_{B_{{\bf i},1} B'_{{\bf i},1}}$ is invertible by Lemma~\ref{gluing-1}.
Therefore we have 
\begin{eqnarray*}
(X^{-1}_{\bf i})_{B'_{{\bf i},1},B_{{\bf i},1}}&=&({\bf 1}_{(\mu-\mu_{0})\times (\mu-\mu_{0})}-(C_{\bf i})_{B_{{\bf i},1},B_{{\bf i},1}})
((I_{\bf i})_{B_{{\bf i},1} B'_{{\bf i},1}})^{-1},\\  
(X^{-1}_{\bf i})_{B'_{{\bf i},1},B_{{\bf i},2}}&=&-((C_{\bf i})_{B_{{\bf i},1},B_{{\bf i},2}})((I_{\bf i})_{B_{{\bf i},1},B'_{{\bf i},1}})^{-1}.  
\end{eqnarray*}
Thus we can determine the representation matrix of $\chi^{-1}$ concerning the $\ZZ$--bases 
$\displaystyle \bigcup_{{\bf i}\in \I} B'_{{\bf i},1}$
and $\displaystyle \bigcup_{{\bf i}\in \I} B_{{\bf i},1}$. 
Therefore, the bilinear form $\chi$ should be unique if it exists.  
\qed
\end{pf}
We have finished the proof of this theorem.
\qed
\end{pf}

\begin{defn}
Let $R=(K_0(R),I_R,\Delta_{re}(R),c_R)$ be a simply-laced generalized root system.
The bilinear form $\chi_R$ in Proposition~\ref{Euler form} is called the {\it Euler form} of $R$. 
\end{defn}

\subsection{Properties of morphisms} In this subsection, we give a sufficient condition for monomorphism. 
Moreover we also give an example of a morphism which is not  monomorphism.

\begin{prop}
Let $\phi: R\rightarrow R'$ be a morphism.
If the Cartan form $I_{R}$ of $R$ is non-degenerate,
then $\phi$ is a monomorphism.
\end{prop}
\begin{pf}
We show that a morphism $\phi$ induces an isometric $\ZZ$--homomorphism for Euler forms:
\begin{lem}
Let $\phi:R\rightarrow R'$ be a morphism.
Assume that the Cartan form $I_{R}$ of $R$ is non-degenerate.
Then we have
\begin{equation*}
\chi_{R'}\left(\phi(\lambda_1),\phi(\lambda_2)\right)=\chi_R(\lambda_1,\lambda_2)
\end{equation*}
for all $\lambda_1,\lambda_2\in K_0(R)$.
\end{lem}
\begin{pf}
Set $\chi'_{R}(\lambda_1,\lambda_2):=\chi_{R'}\left(\phi(\lambda_1),\phi(\lambda_2)\right)$ for all $\lambda_1,\lambda_2\in K_0(R)$.
By the Serre duality for $\chi_{R'}$ and the commutativity $c_{R'}\circ \phi =\phi\circ c_R$, we have 
\begin{eqnarray}
I_{R}(\lambda_1,\lambda_2) &=& I_{R'}\left(\phi(\lambda_1),\phi(\lambda_2)\right) \nonumber\\
&=& \chi_{R'}\left(\phi(\lambda_1),\phi(\lambda_2)\right)+{\chi}_{R'}\left(\phi(\lambda_2),\phi(\lambda_1)\right) \nonumber\\
&=& \chi_{R'}\left(\phi(\lambda_1),\phi(\lambda_2)\right)+{\chi}_{R'}\left(\phi(\lambda_1),-c_{R'}(\phi(\lambda_2))\right)
\nonumber \\
&=& \chi_{R'}\left(\phi(\lambda_1),\phi(\lambda_2)\right)+{\chi}_{R'}\left(\phi(\lambda_1),\phi(-c_{R}(\lambda_2))\right)
\nonumber \\
&=& \chi_{R'}\left(\phi(\lambda_1),\phi((1-c_{R})(\lambda_2))\right)
\nonumber \\
&=& \chi'_{R}(\lambda_1, (1-c_{R})(\lambda_2))
\nonumber \\
&=& \chi'_{R}(\lambda_1, (\chi_{R}^{-1}\circ \chi_{R}\circ (1-c_{R}))(\lambda_2))\label{2.26}\\
&=& \chi'_{R}(\lambda_1, (\chi_{R}^{-1}\circ I_{R})(\lambda_2))\label{2.27}.
\end{eqnarray}
Here we regard the bilinear forms $\chi_{R}$ and $I_{R}$ as elements of ${\rm Hom}_{\ZZ}(K_{0}(R), {\rm Hom}_{\ZZ}(K_{0}(R),\ZZ))$
as in Lemma~\ref{uniqueness-of-euler} in the lines \eqref{2.26} and \eqref{2.27}.
Therefore the non-degeneracy of $I_{R}$ implies $\chi'_{R}\circ \chi_{R}^{-1}={\rm id}_{{\rm Hom}_{\ZZ}(K_{0}(R),\ZZ)}$ and hence $\chi'_{R}=\chi_{R}$.
\qed
\end{pf}
If the $\ZZ$--homomorphism $\phi$ has the nontrivial kernel, it is contradictory to the non-degenerate property of $\chi_{R}$.
Therefore the $\ZZ$--homomorphism $\phi$ is injective. Assume that $\phi\circ \phi_1=\phi\circ \phi_2$ for morphisms $\phi_1, \phi_2:R''\rightarrow R$.
Then the induced $\ZZ$--homomorphism on the lattices $\phi_1, \phi_2: K_{0}(R'')\rightarrow K_{0}(R)$ are equal, and obviously, 
the induced maps on the sets of real roots $\phi_1, \phi_2: \Delta_{re}(R'')\rightarrow \Delta_{re}(R)$ are also equal. 
Therefore the morphism $\phi$ is a monomorphism if $I_{R}$ is non-degenerate.
\qed
\end{pf}

\begin{rem}
If the Cartan form $I_{R}$ is degenerate, there is an example of a morphism which is not a monomorphism.
Let $R$ be the simply-laced generalized root system which is generated by the root basis $(\alpha_{1},\alpha_{2},\alpha_{3})$
satisfying $I_{R}(\alpha_{i},\alpha_{j})=2$ for $i,j=1,2,3$. Moreover we set $R'$ the simply-laced generalized root system 
which is generated by $\beta_1$. Then we can check that
the $\ZZ$-linear map $\phi:K_{0}(R)\rightarrow K_{0}(R')$ mapping $\alpha_{i}$ to $\beta_{1}$ for $i=1,2,3$
induces a morphism which is not a monomorphism by easy calculation.
\end{rem}

\section{Classification of Positive definite generalized root systems}
\subsection{Positive definite generalized root systems}
In this subsection, we shall show irreducible
simply-laced generalized root systems with positive definite Cartan forms 
contain irreducible classical root systems as parts of them. 
\begin{defn}
A simply-laced generalized root system $R=(K_0(R),I_R,\Delta_{re}(R),c_R)$ is called 
{\it positive definite}
if the Cartan form $I_R$ on $K_0(R)\otimes_\ZZ\RR$ is positive definite.
\end{defn} 

\begin{defn}[Section~6.1.1 in \cite{bourbaki}]\label{defn of classical root system}
Let $V$ be a $\mu$-dimensional $\RR$-vector space with an inner product $(-,-)$. 
A subset $\Phi\subset V$ is called a {\it simply-laced classical root system} of $V$ if it satisfies the followings$:$
\begin{enumerate}
\item
$\Phi$ is a finite set such that $(v,v)=2$ for all $v \in \Phi$ and $\RR\Phi=V$.
\item
For $\alpha\in\Phi$, the reflection $r_\alpha$ of $V$ is defined as follows$;$
\begin{equation*}
r_{\alpha}(\lambda) :=\lambda-(\lambda,\alpha)\alpha,\quad \lambda\in V.
\end{equation*}
For $\alpha, \beta\in\Phi$, $r_{\alpha}(\beta)\in\Phi$.
\item
For $\alpha, \beta\in\Phi$, $(\beta,\alpha)\in\ZZ$.
\item
If $\alpha, c\alpha\in\Phi$ for some $c\in\RR$, then $c=\pm1$.
\end{enumerate}
\end{defn}

\begin{prop}
Let $R=(K_0(R),I_R,\Delta_{re}(R),c_R)$ be a positive definite generalized root system.
The set of real roots $\Delta_{re}(R)$ is a classical root system of $K_{0}(R)\otimes_{\ZZ}\RR$ with the inner product $I_R$. 
\end{prop}
\begin{pf}
By the condition {\rm (iii)} in Definition~\ref{object}, 
the condition {\rm (ii)} in Definition~\ref{defn of classical root system} is satisfied.
Since $I_{R}(\alpha, \beta)\in \ZZ$ for all $\alpha, \beta \in \Delta_{re}(R)\subset K_{0}(R)$,
the condition {\rm (ii)} in Definition~\ref{defn of classical root system} is satisfied.
By Lemma~\ref{lem2.4}, the condition {\rm (iv)} in Definition~\ref{defn of classical root system} is satisfied.
By the condition {\rm (i)} and {\rm (ii)} in Definition~\ref{object}, we only need to show the following
\begin{lem}\label{finite}
The set of real roots $\Delta_{re}(R)$ is a finite set.  
\end{lem}
\begin{pf}
Let $B=(\alpha_1, \alpha_2, \dots, \alpha_\mu)$ be a root basis. 
Set $\alpha=c_1\alpha_1+c_2\alpha_2+\cdots +c_\mu\alpha_\mu\in \Delta_{re}(R)\setminus B$
where $c_1,\ c_2\ \dots, c_\mu\in\ZZ$. 
Define the Cartan matrix $C\in M(n,\ZZ)$ as follows$:$
\begin{eqnarray*}
C:=
\begin{pmatrix}
I_{R}(\alpha_1,\alpha_1)&\cdots&I_{R}(\alpha_\mu,\alpha_1)\\
I_{R}(\alpha_1,\alpha_2)&\cdots&I_{R}(\alpha_\mu,\alpha_2)\\
\vdots\\
I_{R}(\alpha_1,\alpha_\mu)&\cdots&I_{R}(\alpha_\mu,\alpha_\mu)\\
\end{pmatrix}.
\end{eqnarray*}
The matrix $C$ is invertible since $I_R$ is positive definite and
$\{\alpha_1, \alpha_2, \dots, \alpha_{\mu}\}$ are linearly independent. 
Then we have
\begin{equation*}
\begin{pmatrix}
c_1\\
c_2\\
\vdots\\
c_\mu\\
\end{pmatrix}=C^{-1}
\begin{pmatrix}
I_{R}(\alpha,\alpha_1)\\
I_{R}(\alpha,\alpha_2)\\
\vdots\\
I_{R}(\alpha,\alpha_\mu)\\
\end{pmatrix}.
\end{equation*}
By the exactly same argument in classical root systems (cf. Section 6.1.3 in \cite{bourbaki}),
we have $-2\leq I_{R}(\alpha, \beta)\leq2$ for all $\alpha, \beta\in\Delta_{re}(R)$, 
in particular, $I_{R}(\alpha, \beta)=\pm2$ if and only if $\alpha=\pm\beta$.
Since $I_{R}(\alpha,\alpha_i)\in\ZZ$, the cardinality of $\Delta_{re}(R)$ is at most $3^{\mu}+\mu$.
\qed
\end{pf}
We have finished the proof of the proposition.
\qed
\end{pf}

\subsection{Admissible diagrams}
In this subsection, we recall admissible diagrams and their properties
introduced and studied by R. W. Carter. Then see \cite{carter} and also \cite{rafael}
for proofs of assertions in this subsection.

We denote by $(-,-)$ an inner product on a $\mu$-dimensional $\RR$-vector space $V$, 
by $\Phi\subset V$ a classical root system and by $W$ the Weyl group
associated to $\Phi$ throughout this and next sections.
\begin{defn}\label{w admi}
Assume that an element $w\in W$ is expressed as 
\begin{equation}\label{3.1.1}
w=(r_{\beta_1} r_{\beta_2} \cdots r_{\beta_{{k_1}}})(r_{\beta_{{k_1}+1}} r_{\beta_{{k_1}+2}} \cdots r_{\beta_{{k_2}}}),
\quad 1\leq {k_1}\leq {k_2}\leq \mu,
\end{equation}
such that
\begin{enumerate}
\item $\beta_1,\dots,\beta_{{k_1}}, \beta_{{k_1}+1},\dots,\beta_{{k_2}}\in \Phi$ are linearly independent over $\RR$,
\item $\beta_1,\dots, \beta_{{k_1}}$ are orthogonal to each other and 
$\beta_{k_1+1},\dots, \beta_{{k_2}}$ are also orthogonal to each other. 
\end{enumerate}
We call an expression \eqref{3.1.1} an {\it admissible representation} of $w$. 
An element $w$ is called {\it admissible} if  an admissible representation of $w$ satisfies
\begin{equation*}
W=\langle r_{\beta_1},\dots,r_{\beta_{{k_1}}}, r_{\beta_{{k_1}+1}},\dots, r_{\beta_{{k_2}}}\rangle.
\end{equation*}
For an admissible element $w\in W$ with an expression \eqref{3.1.1}, 
an {\it admissible diagram} associated to $w$ is defined as follows$:$
\begin{enumerate}
\item the set of vertices is 
$\{\beta_1,\dots,\beta_{{k_1}}$, $\beta_{{k_1}+1},\dots,\beta_{{k_2}}\}$, 
\item the edge between the vertices $\beta_{i}$ and $\beta_{j}$
is given by the following rule$:$
\begin{subequations}
\begin{align*}
\xymatrix{ \circ_{\beta_i}  & \circ_{\beta_j}} & \quad\text{if}\quad (\beta_i,\beta_j)=0,\\
\xymatrix{ \circ_{\beta_i}\ar@{-}[r]  & \circ_{\beta_j}} & \quad\text{if}\quad (\beta_i,\beta_j)\ne 0.
\end{align*}
\end{subequations}
\end{enumerate}
\end{defn}

\begin{prop}[Theorem A in \cite{carter}, Theorem 4.1 and Theorem 6.5 in \cite{rafael}]\label{classifycation}
We have the followings$:$
\begin{enumerate}
\item
Assume that $\Phi$ is irreducible.
For two expressions \eqref{3.1.1} of an admissible element $w\in W$, admissible diagrams
associated to them are isomorphic to each other.
\item
If $\Phi$ is irreducible,
then the admissible diagram associated to $w\in W$ is isomorphic to the one of following diagrams$:$
\begin{eqnarray*}
A_\mu&:&
\begin{xy}
\ar@{-} (0,0) *++!D{} *\cir<3pt>{}="A";
 (10,0) *++!D{} *\cir<3pt>{}="F"
\ar@{-}  "F";
 (20,0) *++!D{} *\cir<3pt>{}="E"
\ar@{-}  "E";(25,0) \ar@{.} (25,0);(35,0)^*!U{}
\ar@{-}  (35,0);
 (40,0) *++!D{} *\cir<3pt>{}="D"
\ar@{-}  "D";
 (50,0) *++!D{} *\cir<3pt>{}
\end{xy} \nonumber \\
D_\mu&:&
\begin{xy}
\ar@{-} (0,0) *++!D{} *\cir<3pt>{};
(10,0) *++!D{} *\cir<3pt>{}="B",
\ar@{-} "B";
(20,0) *++!D{} *\cir<3pt>{}="C"
\ar@{-}"C";(25,0) \ar@{.} (25,0);(35,0)^*!U{}
\ar@{-}  (35,0);
(40,0) *++!L{} *\cir<3pt>{}="D"
\ar@{-} "D";(45,5.55) *++!L{} *\cir<3pt>{}
\ar@{-} "D";(45,-5.55) *++!L{} *\cir<3pt>{}
\end{xy}(4\leq\mu) \nonumber \\
D_{\mu}(a_k)&:&
\begin{xy}
\ar@{-} (0,0) *++!D{k} *\cir<3pt>{}="A";
(5,5.50) *++!{} *\cir<3pt>{}="A_1"
\ar@{-}(0,0) *\cir<3pt>{};
(5,-5.50) *++!{} *\cir<3pt>{}="A_2"
\ar@{-} "A_1";
(10,0) *++!{} *\cir<3pt>{}
\ar@{-} "A_2";
(10,0) *++!{} *\cir<3pt>{}="B"
\ar@{-}"B";(14,0) \ar@{.} (14,0);(19,0)
\ar@{-}  (19,0);
(23,0) *++!D{\mu} *\cir<3pt>{}="D"
\ar@{-} "A";
(-7,0) *++!D{k-1} *\cir<3pt>{}="-B"
\ar@{-}"-B";(-11,0) \ar@{.} (-11,0);(-16,0)
\ar@{-}  (-16,0);
(-20,0) *++!D{2} *\cir<3pt>{}="-D"
\ar@{-} "-D";
(-27,0) *++!D{1} *\cir<3pt>{}
\end{xy} (4\leq\mu, 1\leq k\leq \left[\frac{\mu+1}{2}\right]) \nonumber \\
E_6&:&
\begin{xy}
\ar@{-} (0,0) *++!D{} *\cir<3pt>{}="A";
 (10,0) *++!D{} *\cir<3pt>{}="F"
\ar@{-}  "F";
 (20,0) *++!D{} *\cir<3pt>{}="E"
\ar@{-}  "E";
 (30,0) *++!D{} *\cir<3pt>{}="D"
\ar@{-}  "D"; 
 (40,0) *++!D{} *\cir<3pt>{}="B"
\ar@{-}  "E";
 (20,-10) *++!L{} *\cir<3pt>{}
\end{xy} \nonumber \\
E_{6}(a_1)&:&
\begin{xy}
\ar@{-} (0,0) *++!D{} *\cir<3pt>{}="A";
 (10,0) *++!D{} *\cir<3pt>{}="F"
\ar@{-}  "F";
 (20,0) *++!D{} *\cir<3pt>{}="E"
\ar@{-}  "A";
 (0,-10) *++!D{} *\cir<3pt>{}="D"
\ar@{-}  "D"; 
 (10,-10) *++!D{} *\cir<3pt>{}="B"
\ar@{-} "B";
 (10,0) *++!D{} *\cir<3pt>{}
\ar@{-}  "B";
 (20,-10) *++!L{} *\cir<3pt>{}
\end{xy} \nonumber \\
E_{6}(a_2)&:&
\begin{xy}
\ar@{-} (0,0) *++!D{} *\cir<3pt>{}="A";
 (10,0) *++!D{} *\cir<3pt>{}="F"
\ar@{-}  "F";
 (20,0) *++!D{} *\cir<3pt>{}="E"
\ar@{-}  "A";
 (0,-10) *++!D{} *\cir<3pt>{}="D"
\ar@{-}  "D"; 
 (10,-10) *++!D{} *\cir<3pt>{}="B"
\ar@{-} "B";
 (10,0) *++!D{} *\cir<3pt>{}
\ar@{-}  "B";
 (20,-10) *++!L{} *\cir<3pt>{}="C"
\ar@{-}  "C";
"E"
\end{xy} \nonumber \\
E_7&:&
\begin{xy}
\ar@{-} (0,0) *++!D{} *\cir<3pt>{}="A";
 (10,0) *++!D{} *\cir<3pt>{}="F"
\ar@{-}  "F";
 (20,0) *++!D{} *\cir<3pt>{}="E"
\ar@{-}  "E";
 (30,0) *++!D{} *\cir<3pt>{}="D"
\ar@{-}  "D"; 
 (40,0) *++!D{} *\cir<3pt>{}="B"
 \ar@{-}  "B"; 
 (50,0) *++!D{} *\cir<3pt>{}
\ar@{-}  "E";
 (20,-10) *++!L{} *\cir<3pt>{}
\end{xy} \nonumber \\
E_{7}(a_1)&:&
\begin{xy}
\ar@{-} (0,0) *++!D{} *\cir<3pt>{}="A";
 (10,0) *++!D{} *\cir<3pt>{}="F"
\ar@{-}  "F";
 (20,0) *++!D{} *\cir<3pt>{}="E"
\ar@{-}  "A";
 (0,-10) *++!D{} *\cir<3pt>{}="D"
\ar@{-}  "D"; 
 (10,-10) *++!D{} *\cir<3pt>{}="B"
\ar@{-} "B";
 (10,0) *++!D{} *\cir<3pt>{}
\ar@{-}  "B";
 (20,-10) *++!L{} *\cir<3pt>{}
\ar@{-}  "E";
 (30,0) *++!D{} *\cir<3pt>{}="T" 
\end{xy} \nonumber \\
E_{7}(a_2)&:&
\begin{xy}
\ar@{-} (0,0) *++!D{} *\cir<3pt>{}="A";
 (10,0) *++!D{} *\cir<3pt>{}="F"
\ar@{-}  "F";
 (20,0) *++!D{} *\cir<3pt>{}="E"
\ar@{-}  "A";
 (0,-10) *++!D{} *\cir<3pt>{}="D"
\ar@{-}  "D"; 
 (10,-10) *++!D{} *\cir<3pt>{}="B"
\ar@{-} "B";
 (10,0) *++!D{} *\cir<3pt>{}
\ar@{-}  "B";
 (20,-10) *++!L{} *\cir<3pt>{}
 \ar@{-}  "A";
 (-10,0) *++!D{} *\cir<3pt>{}="V"
\end{xy} \nonumber \\
E_{7}(a_3)&:&
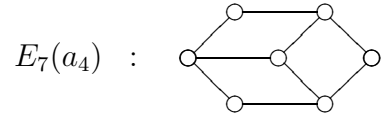
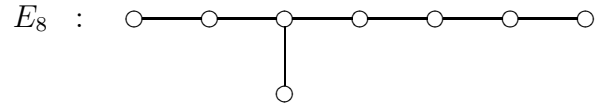
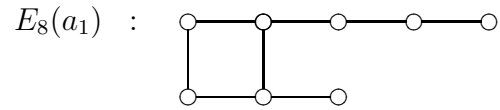
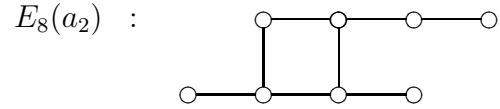
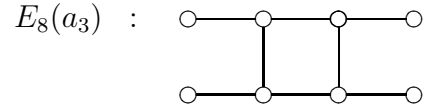
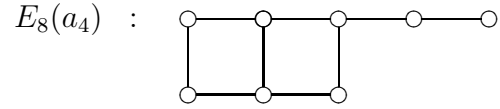
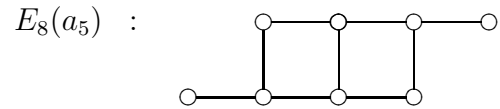
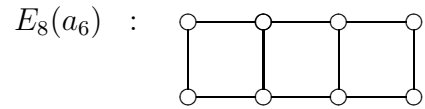
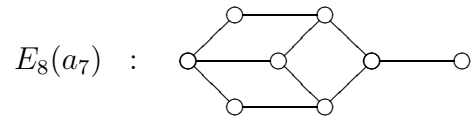
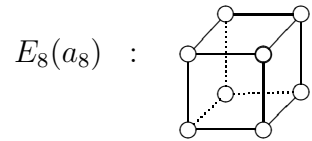
\begin{xy}
\ar@{-} (0,0) *++!D{} *\cir<3pt>{}="A";
 (10,0) *++!D{} *\cir<3pt>{}="F"
\ar@{-}  "F";
 (20,0) *++!D{} *\cir<3pt>{}="E"
\ar@{-}  "A";
 (0,-10) *++!D{} *\cir<3pt>{}="D"
\ar@{-}  "D"; 
 (10,-10) *++!D{} *\cir<3pt>{}="B"
\ar@{-} "B";
 (10,0) *++!D{} *\cir<3pt>{}
\ar@{-}  "B";
 (20,-10) *++!L{} *\cir<3pt>{}="C"
\ar@{-}  "C";
"E"
\ar@{-}  "E";
 (30,0) *++!D{} *\cir<3pt>{}="W"
\end{xy} \nonumber \\
\end{eqnarray*}
\begin{eqnarray*}
E_{7}(a_4)&:&
\begin{xy}
\ar@{-} (0,0) *++!{} *\cir<3pt>{}="A";
(6.2,6.100) *++!{} *\cir<3pt>{}="A_1"
\ar@{-}(0,0) *\cir<3pt>{};
(6.2,-6.100) *++!{} *\cir<3pt>{}="A_2"
\ar@{-} "A";
(12,0) *++!{} *\cir<3pt>{}="B"
\ar@{-}  "A_1";
(18.2,6.100) *++!{} *\cir<3pt>{}="C"
\ar@{-}  "A_2";
(18.2,-6.100) *++!{} *\cir<3pt>{}="D"
\ar@{-}  "B";
"C"
\ar@{-}  "B";
"D"
\ar@{-}  "C";
(24.4,0) *++!{} *\cir<3pt>{}
\ar@{-}  "D";
(24.4,0) *++!{} *\cir<3pt>{}
\end{xy} \nonumber \\
E_8&:&
\begin{xy}
\ar@{-} (0,0) *++!D{} *\cir<3pt>{}="A";
 (10,0) *++!D{} *\cir<3pt>{}="F"
\ar@{-}  "F";
 (20,0) *++!D{} *\cir<3pt>{}="E"
\ar@{-}  "E";
 (30,0) *++!D{} *\cir<3pt>{}="D"
\ar@{-}  "D"; 
 (40,0) *++!D{} *\cir<3pt>{}="B"
\ar@{-} "B";
 (50,0) *++!D{} *\cir<3pt>{}="C"
\ar@{-} "C";
 (60,0) *++!D{} *\cir<3pt>{}
\ar@{-}  "E";
 (20,-10) *++!L{} *\cir<3pt>{}
\end{xy} \nonumber \\
E_{8}(a_1)&:&
\begin{xy}
\ar@{-} (0,0) *++!D{} *\cir<3pt>{}="A";
 (10,0) *++!D{} *\cir<3pt>{}="F"
\ar@{-}  "F";
 (20,0) *++!D{} *\cir<3pt>{}="E"
\ar@{-}  "A";
 (0,-10) *++!D{} *\cir<3pt>{}="D"
\ar@{-}  "D"; 
 (10,-10) *++!D{} *\cir<3pt>{}="B"
\ar@{-} "B";
 (10,0) *++!D{} *\cir<3pt>{}
\ar@{-}  "B";
 (20,-10) *++!L{} *\cir<3pt>{}
\ar@{-}  "E";
 (30,0) *++!D{} *\cir<3pt>{}="T"
 \ar@{-}  "T";
 (40,0) *++!D{} *\cir<3pt>{}
\end{xy} \nonumber \\
E_{8}(a_2)&:&
\begin{xy}
\ar@{-} (0,0) *++!D{} *\cir<3pt>{}="A";
 (10,0) *++!D{} *\cir<3pt>{}="F"
\ar@{-}  "F";
 (20,0) *++!D{} *\cir<3pt>{}="E"
\ar@{-}  "A";
 (0,-10) *++!D{} *\cir<3pt>{}="D"
\ar@{-}  "D"; 
 (10,-10) *++!D{} *\cir<3pt>{}="B"
\ar@{-} "B";
 (10,0) *++!D{} *\cir<3pt>{}
\ar@{-}  "B";
 (20,-10) *++!L{} *\cir<3pt>{}
\ar@{-}  "D";
(-10,-10) *++!D{} *\cir<3pt>{}
\ar@{-}  "E";
 (30,0) *++!D{} *\cir<3pt>{}
\end{xy} \nonumber \\
E_{8}(a_3)&:&
\begin{xy}
\ar@{-} (0,0) *++!D{} *\cir<3pt>{}="A";
 (10,0) *++!D{} *\cir<3pt>{}="F"
\ar@{-}  "F";
 (20,0) *++!D{} *\cir<3pt>{}="E"
\ar@{-}  "A";
 (0,-10) *++!D{} *\cir<3pt>{}="D"
\ar@{-}  "D"; 
 (10,-10) *++!D{} *\cir<3pt>{}="B"
\ar@{-} "B";
 (10,0) *++!D{} *\cir<3pt>{}
\ar@{-}  "B";
 (20,-10) *++!L{} *\cir<3pt>{}
\ar@{-}  "A";
 (-10,0) *++!D{} *\cir<3pt>{}="V"
\ar@{-}  "D";
 (-10,-10) *++!D{} *\cir<3pt>{}
\end{xy} \nonumber \\
E_{8}(a_4)&:&
\begin{xy}
\ar@{-} (0,0) *++!D{} *\cir<3pt>{}="A";
 (10,0) *++!D{} *\cir<3pt>{}="F"
\ar@{-}  "F";
 (20,0) *++!D{} *\cir<3pt>{}="E"
\ar@{-}  "A";
 (0,-10) *++!D{} *\cir<3pt>{}="D"
\ar@{-}  "D"; 
 (10,-10) *++!D{} *\cir<3pt>{}="B"
\ar@{-} "B";
 (10,0) *++!D{} *\cir<3pt>{}
\ar@{-}  "B";
 (20,-10) *++!L{} *\cir<3pt>{}="C"
\ar@{-}  "C";
"E"
\ar@{-}  "E";
 (30,0) *++!D{} *\cir<3pt>{}="W"
 \ar@{-}  "W";
 (40,0) *++!D{} *\cir<3pt>{}
\end{xy} \nonumber \\
E_{8}(a_5)&:&
\begin{xy}
\ar@{-} (0,0) *++!D{} *\cir<3pt>{}="A";
 (10,0) *++!D{} *\cir<3pt>{}="F"
\ar@{-}  "F";
 (20,0) *++!D{} *\cir<3pt>{}="E"
\ar@{-}  "A";
 (0,-10) *++!D{} *\cir<3pt>{}="D"
\ar@{-}  "D"; 
 (10,-10) *++!D{} *\cir<3pt>{}="B"
\ar@{-} "B";
 (10,0) *++!D{} *\cir<3pt>{}
\ar@{-}  "B";
 (20,-10) *++!L{} *\cir<3pt>{}="C"
\ar@{-}  "C";
"E"
\ar@{-}  "E";
 (30,0) *++!D{} *\cir<3pt>{}
 \ar@{-}  "D";
 (-10,-10) *++!D{} *\cir<3pt>{}
\end{xy} \nonumber \\
E_{8}(a_6)&:&
\begin{xy}
\ar@{-} (0,0) *++!D{} *\cir<3pt>{}="A";
 (10,0) *++!D{} *\cir<3pt>{}="F"
\ar@{-}  "F";
 (20,0) *++!D{} *\cir<3pt>{}="E"
\ar@{-}  "A";
 (0,-10) *++!D{} *\cir<3pt>{}="D"
\ar@{-}  "D"; 
 (10,-10) *++!D{} *\cir<3pt>{}="B"
\ar@{-} "B";
 (10,0) *++!D{} *\cir<3pt>{}
\ar@{-}  "B";
 (20,-10) *++!L{} *\cir<3pt>{}="C"
\ar@{-}  "C";
"E"
\ar@{-}  "E";
 (30,0) *++!D{} *\cir<3pt>{}="W"
\ar@{-}  "W";
 (30,-10) *++!D{} *\cir<3pt>{}="T"
\ar@{-} "T";
"C"
\end{xy} \nonumber \\
E_{8}(a_7)&:&
\begin{xy}
\ar@{-} (0,0) *++!{} *\cir<3pt>{}="A";
(6.2,6.100) *++!{} *\cir<3pt>{}="A_1"
\ar@{-}(0,0) *\cir<3pt>{};
(6.2,-6.100) *++!{} *\cir<3pt>{}="A_2"
\ar@{-} "A";
(12,0) *++!{} *\cir<3pt>{}="B"
\ar@{-}  "A_1";
(18.2,6.100) *++!{} *\cir<3pt>{}="C"
\ar@{-}  "A_2";
(18.2,-6.100) *++!{} *\cir<3pt>{}="D"
\ar@{-}  "B";
"C"
\ar@{-}  "B";
"D"
\ar@{-}  "C";
(24.4,0) *++!{} *\cir<3pt>{}
\ar@{-}  "D";
(24.4,0) *++!{} *\cir<3pt>{}="E"
\ar@{-} "E";
(36.4,0) *++!{} *\cir<3pt>{}
\end{xy} \nonumber \\
E_{8}(a_8)&:&
\begin{xy}
\ar@{-} (0,0) *++!D{} *\cir<3pt>{}="A";
 (10,0) *++!D{} *\cir<3pt>{}="F"
\ar@{-}  "A";
 (0,-10) *++!D{} *\cir<3pt>{}="D"
\ar@{-}  "D"; 
 (10,-10) *++!D{} *\cir<3pt>{}="B"
\ar@{-} "B";
 (10,0) *++!D{} *\cir<3pt>{}
\ar@{-} "A";
 (5,5.30) *++!{} *\cir<3pt>{}="A_1"
\ar@{-} "A_1";
 (15,5.30) *++!{} *\cir<3pt>{}="F_1"
\ar@{-} "F_1";
 (15,-5) *++!D{} *\cir<3pt>{}="F_2"
\ar@{.} "D";
 (5,-5.30) *++!{} *\cir<3pt>{}="D_1"
\ar@{-} "B";
"F_2"
\ar@{-} "F";
"F_1"
\ar@{.} "A_1";
"D_1"
\ar@{.} "F_2";
"D_1"
\end{xy} \nonumber \\
\end{eqnarray*}
\end{enumerate}
\end{prop}

The next proposition follows from Theorem A, Theorem B, Proposition 21 in \cite{carter}$:$
\begin{prop}\label{admi uniqueness}
We have the followings$:$
\begin{enumerate}
\item For a diagram $A$ in Proposition~\ref{classifycation} {\rm (ii)}, 
there exists a pair of a classical root system and an admissible element $(\Phi, w)$ such that 
the admissible diagram associated to $w$ is isomorphic to $A$.
\item If a pair $(\Phi,w)$ gives the diagram $A_{\mu}$, $D_{\mu}$ or $D_{\mu}(a_{k})$, 
$E_{6}$ or $E_{6}(a_{k})$, $E_{7}$ or $E_{7}(a_{k})$, $E_{8}$ or $E_{8}(a_{k})$ 
in Proposition~\ref{classifycation}, then $\Phi$ is the one of types $A_{\mu}$
$D_{\mu}$, $E_{6}, E_{7}, E_{8}$ respectively.
\item If the admissible diagram associated to $w$ is isomorphic to the one associated to $w'$ for
pairs $(\Phi, w)$ and $(\Phi', w')$, then $\Phi$ is isomorphic to $\Phi'$. 
\end{enumerate}
\end{prop}

\subsection{Classification}
Admissible diagrams by R. W. Carter associated to their Coxteter elements
enable us to classify positive definite generalized root systems. 

We denote by $R=(K_0(R),I_R,\Delta_{re}(R),c_R)$
an irreducible positive definite generalized root systems
throughout this subsection.

The next proposition follows from Lemma 3, Lemma 5 and Theorem C in \cite{carter}$:$ 
\begin{prop}\label{Carter cor}
Any element in $W$ can be expressed by an admissible representation.
\end{prop}

The following is the key proposition to connect irreducible positive definite generalized root systems
with admissible diagrams.
\begin{prop}\label{admissible}
The Coxeter element $c_{R}$ is admissible.
\end{prop}
\begin{pf}
Take an admissible representation  
$c_{R}=(r_{\beta_1} r_{\beta_2} \cdots r_{\beta_{k}})(r_{\beta_{k+1}} r_{\beta_{k+2}} \cdots r_{\beta_{k'}})$
of the Coxeter element. 

\begin{lem}[Lemma 3 in \cite{carter}]\label{lem3 in carter}
For an element $w'\in W$,
set ${\rm l}(w'):={\rm min}\{ k\in \ZZ ~\vert~ w'=
r_{\alpha_1}r_{\alpha_2}\cdots r_{\alpha_k},\ \alpha_1, \alpha_2, \dots, \alpha_k \in \Phi \}$.
For the element $w=r_{\gamma_{1}}\cdots r_{\gamma_{k}}\in W$, one has $k={\rm l}(w)$ if and only if
$\gamma_{1},\dots,\gamma_{k}\in \Phi$ are linearly independent.
\end{lem}
Denote a root basis of $R$ by $(\alpha_1, \dots , \alpha_{\mu})$. 
It is obvious that $k'\le \mu$ since $\beta_{1},\dots, \beta_{k'}$ are linearly independent and $(\alpha_{1},\dots,\alpha_{\mu})$
are root basis.
If $k'< \mu$, the elements $\alpha_{1},\dots \alpha_{\mu}$ would be linearly dependent by Lemma~\ref{lem3 in carter}.
However the elements $\alpha_{1},\dots \alpha_{\mu}$ must be linear independent since they form a root basis.
Therefore an admissible representation of $c_{R}$ is written as 
$c_R =(r_{\beta_1} r_{\beta_2} \cdots r_{\beta_{k}})(r_{\beta_{k+1}} r_{\beta_{k+2}} \cdots r_{\beta_{\mu}})$
where $1\leq k\leq \mu$.

\begin{lem}\label{realization of admi}
Take an element $c\in W$ and an admissible representation of $c$,
$c=(r_{\beta_1}\cdots r_{\beta_{{k_1}}})(r_{\beta_{{k_1}+1}}\cdots r_{\beta_{{k_2}}})$
where $1\leq {k_1}\leq {k_2}\leq \mu$. 
There exists the positive definite generalized root system $R'$ whose root basis is 
the ordered set $(\beta_1,\dots,\beta_{{k_1}}$, $\beta_{{k_1}+1},\dots,\beta_{{k_2}})$ 
and, in particular, whose Coxeter element is $c$.   
\end{lem}
\begin{pf}
Define the simply-laced generalized root system $R'=({K_0}(R'),I_{R'},\Delta_{re}(R'),c_{R'})$ as follows$:$
\begin{eqnarray*}
{K_0}(R')&:=&\ZZ\beta_1+\dots +\ZZ\beta_{{k_2}}, \\
I_{R'} (\beta_i, \beta_j)&:=&I_{R}(\beta_i, \beta_j)\quad  (1\leq i\neq i\leq {k_2}), \\
W(R')&:=&\langle r_{\beta_1},\dots,r_{\beta_{{k_2}}}\rangle, \\
\Delta_{re}(R')&:=&W(R')\{\beta_1,\dots,\beta_{{k_2}}\}, \\
c_{R'}&:=&(r_{\beta_1} r_{\beta_2} \cdots r_{\beta_{{k_1}}}) (r_{\beta_{{k_1}+1}} r_{\beta_{{k_1}+2}} \cdots r_{\beta_{{k_2}}}).
\end{eqnarray*}
The elements $\beta_1,\dots,\beta_{{k_1}}$, $\beta_{{k_1}+1},\dots,\beta_{{k_2}}$
are linear independent over $\RR$ by Definition~\ref{w admi}. 
It is obvious that a root basis of $R'$ is $(\beta_1,\dots,\beta_{{k_1}}$, $\beta_{{k_1}+1},\dots,\beta_{{k_2}})$. 
Moreover $R'$ is positive definite since $R$ is so.
It is also obvious that $c_{R'}$ is the Coxeter element of $R'$.
\qed
\end{pf}

Lemma~\ref{realization of admi} gives 
the natural simply-laced generalized root system $R'$ whose root basis is 
$(\beta_1,\dots,\beta_{k}$, $\beta_{k+1},\dots,\beta_{\mu})$.
Then we have the natural monomorphism $\iota:R'\longrightarrow R$
since the root basis $(\beta_1,\dots,\beta_{k}$, $\beta_{k+1},\dots,\beta_{\mu})$ is a subset of $\Delta_{re}(R)$ and $c_R=c_{R'}$. 

Let ${\chi_R}|_{R'}$ be the restriction of $\chi_{R}$ to $K_{0}(R')$. By Theorem~\ref{Euler form}, we have
\begin{equation*}
{\chi_R}|_{R'}=\chi_{R'}.
\end{equation*}
Let $A, B\in M_{\mu}(\RR)$ be the matrix representations of ${\chi_R}|_{R'}, \chi_{R'}$ under the root beses above
respectively.
Since $(\beta_1,\dots,\beta_{\mu})\subset K_{0}(R)$, there exists the matrix $P$ satisfying 
\begin{equation*}
(\beta_1,\dots,\beta_{\mu})=(\alpha_1,\dots , \alpha_{\mu})P, \quad P\in M_{\mu}(\ZZ),
\end{equation*}
\begin{equation*}
B=P^{T}AP.
\end{equation*}
Then $P\in {\rm GL}({\mu},\ZZ)$ since $\det{A}=\det{B}=1$ and $\det{P}=\pm1$.
We have $K_0(R)=K_0(R')$ and $I_{R}=I_{R'}$. Hence $R$ and $R'$ have the sets of real roots which is
isomorphic to the classical root system defined by $I_{R}$. 
Therefore we have $W(R)=W(R')$ and $c_R$ is admissible.
\qed
\end{pf}

By Proposition~\ref{admissible}, the Coxeter element is admissible.
Hence we can choose
the diagram $A_{R}$ in Proposition~\ref{classifycation} {\rm (ii)} 
which is isomorphic to the admissible diagram associated to the Coxeter element $c_{R}$ of $R$.
Therefore we have the map $\varphi$ from
the set $\R^{{\rm irr, pd}}$ of irreducible positive definite generalized root systems to
the set $\Gamma$ of diagrams in Proposition~\ref{classifycation} {\rm (ii)}$:$
\begin{equation*}
\varphi:\R^{{\rm irr, pd}}\longrightarrow \Gamma, \quad R\mapsto A_{R}.
\end{equation*}

We denote the natural equivalence relation on $\R^{{\rm irr, pd}}$ by $R\sim R'$ if and only if $R$ is isomorphic to $R'$.
Under this notation, we have the main theorem of this section$:$
\begin{thm}\label{classification thm}
The map
\begin{equation*}
\varphi:\R^{{\rm irr, pd}}\longrightarrow \Gamma, \quad R\mapsto A_{R},
\end{equation*}
induces the  natural bijection$:$
\begin{equation*}
\overline{\varphi}:\R^{{\rm irr, pd}}/\sim \overset{\cong}\longrightarrow \Gamma, \quad [R]\mapsto A_{R},
\end{equation*}
where $[R]$ is the isomorphism class of $R \in \R^{{\rm irr, pd}}$.
\end{thm}
\begin{pf}
We divide the proof of this theorem into the following three steps.

\noindent
\underline{\bf {\rm (i)} $\overline{\varphi}$ is well-defined.}
\vspace{5pt}

\begin{lem}\label{a}
If $R$ is isomorphic to $R'$, 
then the admissible diagram associated to $c_{R}$ is isomorphic to the one associated to $c_{R'}$.
\end{lem}
\begin{pf}
Assume that there exists an isomorphism $\phi:R\longrightarrow R'$. 
We denote by $c_R =(r_{\alpha_1} r_{\alpha_2} \cdots r_{\alpha_{k}}) (r_{\alpha_{k+1}} \cdots r_{\alpha_{\mu}})$
an admissible representation of $c_R$.
By Lemma~\ref{commute c and phi}, we have 
\begin{eqnarray}\label{0.5.2}
\phi \circ c_R \circ \phi^{-1}&=&\phi \circ (r_{\alpha_1} r_{\alpha_2} \cdots r_{\alpha_{k}} r_{\alpha_{k+1}} \cdots r_{\alpha_{\mu}})
 \circ \phi^{-1} \nonumber \\ 
&=&(\phi \circ r_{\alpha_1} \circ \phi^{-1} )(\phi \circ r_{\alpha_2} \circ \phi^{-1})\cdots 
(\phi \circ r_{\alpha_{\mu}} \circ \phi^{-1}) \nonumber \\ 
&=& r_{\phi(\alpha_1)}r_{\phi(\alpha_2)}\cdots r_{\phi(\alpha_{\mu})} \nonumber \\ 
&=& (r_{\phi(\alpha_1)}r_{\phi(\alpha_2)}\cdots r_{\phi(\alpha_{k})}) ( r_{\phi(\alpha_{k+1})} \cdots r_{\phi(\alpha_{\mu})}).
\end{eqnarray}
By Definition~\ref{morphism} {\rm (iii)}, we have 
$c_{R'}=(r_{\phi(\alpha_1)}r_{\phi(\alpha_2)}\cdots r_{\phi(\alpha_{k})}) ( r_{\phi(\alpha_{k+1})} \cdots r_{\phi(\alpha_{\mu})})$.
Moreover the representation \eqref{0.5.2} should be an admissible representation of $c_{R'}$ by Definition~\ref{morphism} {\rm(i)}. 
Again by Definition~\ref{morphism} {\rm(i)}, the admissible diagram associated to $c_R$ is isomorphic to  
the one associated to $c_{R'}$ under above admissible representations.
\qed
\end{pf}

\noindent
\underline{\bf {\rm (ii)} $\overline{\varphi}$ is injective.}
\vspace{5pt}

The following lemma is necessary to show that $\overline{\varphi}$ is injective. 
\begin{lem}\label{b}
If the admissible diagram associated to $c_{R}$ is isomorphic to the admissible diagram associated to $c_{R'}$,
then there exists an $\RR$--isometric isomorphism $f:K_0(R)\otimes \RR\longrightarrow K_0(R')\otimes \RR$ 
such that $f(\Delta_{re}(R))=\Delta_{re}(R')$. 
\end{lem}
\begin{pf}
By Proposition~\ref{admi uniqueness}, the classical root system obtained from 
the admissible diagram associated to $c_{R}$ is isomorphic to the one obtained from 
the admissible diagram associated to $c_{R'}$.
Then there exists an isomorphism between $\Delta_{re}(R)$ and $\Delta_{re}(R')$.
Denote by $\{ \alpha_1,\dots , \alpha_{\mu}\}$ a classical root basis of 
the classical root system obtained from $R$.
Here we can choose a classical root basis
$\{ \beta_1,\dots,\beta_{\mu} \}$ of the classical root system obtained from $R'$ such that
\begin{equation*}
I_{R'}(\beta_i,\beta_j)=I_{R}(\alpha_i, \alpha_j) \quad (1\leq i\neq j\leq \mu).
\end{equation*}
This Lemma holds.
\qed
\end{pf}

\begin{prop}[Theorem 6.5 in \cite{rafael}]\label{admi conjugate}
Assume that $\Phi$ is irreducible and $w_1$, $w_2\in W$ are admissible
elements whose admissible diagrams are isomorphic to each other. 
If $\Phi$ is the one of types $A_{\mu}, D_{\mu} , E_6, E_7, E_8$, then $w_1$ and $w_2$ are conjugate. 
\end{prop}

Combining above lemma and proposition, we show that $\overline{\varphi}$ is injective.
\begin{lem}\label{a}
If the admissible diagram associated to $c_{R}$ is isomorphic to the one associated to $c_{R'}$,
then $R$ is isomorphic to $R'$. 
\end{lem}
\begin{pf}
Assume that the admissible diagram associated to $c_R$ is isomorphic to the one associated to $c_{R'}$.
We denote the admissible representations of $c_R$ and $c_{R'}$ by
$c_R=(r_{\alpha_1} r_{\alpha_2} \cdots r_{\alpha_{k}}) (r_{\alpha_{k+1}} \cdots r_{\alpha_{\mu}})$ and
$c_{R'}=(r_{\beta_1} r_{\beta_2} \cdots r_{\beta_{k}}) (r_{\beta_{k+1}} \cdots r_{\beta_{\mu}})$ respectively.

By Lemma~\ref{b},
there exists an element $c \in W(R)$, whose admissible representation is given by
$(r_{\beta_1'} r_{\beta_2'} \cdots r_{\beta_{k}'}) (r_{\beta_{k+1}'} \cdots r_{\beta_{\mu}'})$, satisfying the following;
\begin{equation*}
I_{R}(\beta'_i,\beta'_j)=I_{R'}(\beta_i,\beta_j)\quad (1\leq i\neq j\leq \mu).
\end{equation*}
Lemma~\ref{realization of admi} gives the natural generalized root system $R''$ 
whose root basis is $(\beta_1',\dots,\beta_{\mu}')$. 
Then $R''$ is isomorphic to $R'$. Indeed, if we set a linear map $\phi$ by
\begin{equation*}
\phi:K_0(R'')\longrightarrow K_0(R'),\quad \beta_i' \longmapsto \beta_i,
\end{equation*}
it is obvious from the construction of $(\beta_1',\dots,\beta_{\mu}')$ that $\phi$ defines a morphism
and $\phi$ is isomorphism. 
Here the admissible diagram associated to $c_R$ is isomorphic to the one associated to $c$.
By Proposition~\ref{admi conjugate}, $c_R$ and $c$ are conjugate in $W(R)$, namely, there exists an element
$w\in W(R)$ satisfying 
\begin{equation*}
c_R=w^{-1}cw.
\end{equation*}
If we define the $\ZZ$-linear map as follows;
\begin{equation*}
\psi:K_0(R)\longrightarrow K_0(R'),\quad \alpha \longmapsto w(\alpha),
\end{equation*}
then $\psi$ defines a morphism and $\psi$ is the isomorphism.
Therefore $\psi \circ \phi:R\longrightarrow R'$ is the isomorphism between $R$ and $R'$.
\qed
\end{pf}

\noindent
\underline{\bf {\rm (iii)} $\overline{\varphi}$ is surjective.}
\vspace{5pt}

\begin{lem}
For a diagram $A$ in Proposition~\ref{classifycation} {\rm (ii)}, 
there is an irreducible positive definite generalized root system $R$ such that an admissible diagram
associated to $c_{R}$ is isomorphic to $A$.
\end{lem}
\begin{pf}
This lemma follows from Proposition~\ref{admi uniqueness} and Lemma~\ref{realization of admi}.
\qed
\end{pf}

We have finished the proof of this theorem.
\qed
\end{pf}


\end{document}